\documentclass[11pt,oneside,reqno]{amsart}
\usepackage[all]{xy}
\usepackage{amsmath,amsthm,amsfonts,amssymb,times}
\usepackage[mathscr]{eucal}
\usepackage{verbatim}
\usepackage{url}
\setbox0=\hbox{$+$}
\newdimen\plusheight
\plusheight=\ht0
\def\+{\;\lower\plusheight\hbox{$+$}\;}

\setbox0=\hbox{$-$}
\newdimen\minusheight
\minusheight=\ht0
\def\-{\;\lower\minusheight\hbox{$-$}\;}

\setbox0=\hbox{$\cdots$}
\newdimen\cdotsheight
\cdotsheight=\plusheight
\def\cds{\lower\cdotsheight\hbox{$\cdots$}}

\makeatletter
\def\leqalignno#1{\displ@y \tabskip\z@ plus\@ne fil
  \halign to\displaywidth{\hfil$\@lign\displaystyle{##}$\tabskip\z@skip
    &$\@lign\displaystyle{{}##}$\hfil\tabskip\z@ plus\@ne fil
    &\kern-\displaywidth\rlap{$\@lign\hbox{\rm##}$}\tabskip\displaywidth\crcr
    #1\crcr}}
\makeatother
\let\dotlessi=\i

\newcommand{\eb}{\begin{equation}}
\newcommand{\ee}{\end{equation}}
\renewcommand{\Im}{\operatorname{Im}}
\newcommand{\df}{\dfrac}
\newcommand{\tf}{\tfrac}

\renewcommand{\Re}{\operatorname{Re}}
\renewcommand{\Im}{\operatorname{Im}}

 \renewcommand{\a}{\alpha}
\renewcommand{\b}{\beta}

\newcommand{\G}{\Gamma}
\renewcommand{\l}{\lambda}

\renewcommand{\t}{\tau}

\renewcommand{\Re}{\textup{Re}}
\renewcommand{\Im}{\textup{Im}}

\renewcommand{\(}{\left\(}
\renewcommand{\)}{\right\)}
\renewcommand{\[}{\left\[}
\renewcommand{\]}{\right\]}

\numberwithin{equation}{section}
 \theoremstyle{plain}
\newtheorem{theorem}{Theorem}[section]

\allowdisplaybreaks

\begin{document}
\title[Modified Bessel Functions]
{Modified Bessel Functions in Analytic Number Theory}
\author{Bruce C.~Berndt, Atul Dixit, Rajat Gupta, Alexandru Zaharescu}
\address{Department of Mathematics, University of Illinois, 1409 West Green
Street, Urbana, IL 61801, USA} \email{berndt@illinois.edu}
\address{Department of Mathematics, Indian Institute of Technology Gandhinagar, Palaj, Gandhinagar 382355, Gujarat, India}\email{adixit@iitgn.ac.in}
\address{Institute of Mathematics, Academia Sinica, 6F, Astronomy-Mathematics Building, No. 1, Sec. 4, Roosevelt Road, Taipei 106319, TAIWAN}\email{rajatgpt@gate.sinica.edu.tw}
\address{Department of Mathematics, University of Illinois, 1409 West Green
Street, Urbana, IL 61801, USA; Institute of Mathematics of the Romanian
Academy, P.O.~Box 1-764, Bucharest RO-70700, Romania}
\email{zaharesc@illinois.edu}

\dedicatory{In Memory of Richard A.~Askey}

\begin{abstract}
The modified Bessel functions $K_{\nu}(z)$, or, for brevity, K-Bessel functions, arise at key places in analytic number theory.  In particular, they  appear  in beautiful arithmetic identities.  A survey of these arithmetical identities and their appearances in number theory is provided.
\end{abstract}
\maketitle

\section{Reflections of the First Author}

 I became acquainted with Richard Askey at the University of Wisconsin while writing my Ph.D.~thesis \cite{phdthesis}, wherein Bessel functions play a prominent role.  Our friendship began either when I sought help from him about Bessel functions, or when (more probably) he learned that I was heavily involved with Bessel functions and wanted to see what I was doing with them.  Although I cannot pinpoint any particular facts about Bessel functions that I learned from Dick, what I learned was substantial, and he certainly pointed me to Watson's treatise \cite{watson} to learn more.  In summary, our friendship began in 1965 during the writing of my doctoral dissertation.

 In five of the six chapters of my thesis, the ordinary Bessel functions $J_{\nu}(x)$ are featured, while in Chapter 6, the modified Bessel functions $K_{\nu}(x)$ are the focus.  The first paper that I published appeared in the \emph{Proceedings of the Edinburgh Mathematical Society}, and it was based on Theorem 6.1 in Chapter 6 in my thesis.  In the following survey describing those instances in analytic number theory or related analysis wherein the Bessel functions  $K_{\nu}(x)$ appear, Theorem \ref{thm3} is a generalization of Theorem 6.1 in my thesis \cite{phdthesis}.

\section{Introduction} The ordinary Bessel function $J_{\nu}(z)$, the Bessel function $Y_{\nu}(z)$ of the second kind, the Bessel function of imaginary argument $I_{\nu}(z)$, and the modified Bessel function $K_{\nu}(z)$ are defined, respectively, by \cite[pp.~40, 64, 77, 78]{watson}
\begin{align}
J_{\nu}(z)&:=\sum_{n=0}^{\infty}\df{(-1)^n\left(\frac12 z\right)^{\nu+2n}}{n!\Gamma(\nu+n+1)}, \quad z\in\mathbb{C}, \label{bb}\\
\label{b.Y}
		Y_{\nu}(z)&:=\df{J_{\nu}(z)(\nu\pi)-J_{-\nu}(z)}{\sin(\nu\pi)}, \quad z\in\mathbb{C}, \nu\notin\mathbb{Z},\\
\label{bessely-int}
	Y_{n}(z)&:=\lim_{\nu\to n}Y_{\nu}(z), \quad n\in \mathbb{Z},\\
I_{\nu}(z)&:=\sum_{n=0}^{\infty}\df{(\frac12 z)^{\nu+2n}}{n!\Gamma(\nu+n+1)},\quad z\in \mathbb{C},\label{b1}\\
K_{\nu}(z)&:=\frac{\pi}{2}\df{I_{-\nu}(z)-I_{\nu}(z)}{\sin \nu\pi}, \quad z\in \mathbb{C}, \nu\notin\mathbb{Z},\label{b2}\\
K_n(z)&:=\lim_{\nu\to n}K_{\nu}(z), \quad n\in\mathbb{Z}.\label{b66}
\end{align}
As $z\to\infty$, $J_{\nu}(z)$, $Y_{\nu}(z)$, $I_{\nu}(z)$, and $K_{\nu}(z)$ possess the asymptotic formulas \cite[pp.~199, 202]{watson},
\begin{align}
 J_{\nu}(z)=&\sqrt{\df{2}{\pi z}}\left(\cos(z-\tfrac12 \nu\pi-\tfrac14 \pi)+O\left(\df{1}{z}\right)\right),\label{b5a}\\
 Y_{\nu}(z)=&\sqrt{\df{2}{\pi z}}\left(\sin(z-\tfrac12 \nu\pi-\tfrac14 \pi)+O\left(\df{1}{z}\right)\right),\label{b5c}\\
 I_{\nu}(z)=&\sqrt{\df{1}{2\pi z}}e^z\left(1+O\left(\df{1}{z}\right)\right),\label{b5}\\
K_{\nu}(z)=&\sqrt{\df{\pi}{2 z}}e^{-z}\left(1+O\left(\df{1}{z}\right)\right).\label{b6}
\end{align}
Both $J_{\nu}(z)$ and $K_{\nu}(z)$ arise in several problems and identities in analytic number theory.  In particular, as evinced by the cosine-function appearing in \eqref{b5a}, ordinary Bessel functions appear in certain problems involving circular, spherical, or cylindrical symmetry.  The rapid decay of modified Bessel functions, as indicated in \eqref{b6}, point to possible appearances in summation formulas and contexts where theta functions might otherwise appear.  As our title indicates, in this paper we devote attention to problems and identities in number theory wherein modified Bessel functions appear.

Throughout this paper, $\sigma=\Re(s)$.

\section{The Dirichlet Divisor Problem}\label{ddp}  Let $d(n)$ denote the number of positive divisors of the positive integer $n$. Set
\begin{equation*}\label{d(n)}
D(x):={\sum_{n\leq x}}^\prime d(n),
\end{equation*}
where the prime $\prime$ on the summation sign indicates that if $x$ is an integer, then only $\frac12 d(x)$ is counted.
Historically, perhaps the first instance in number theory in which modified Bessel functions appear is in the ``Vorono\"{\dotlessi} Identity,'' first proved by  G.~F.~Vorono\"{\dotlessi} \cite{voronoi} in 1904.  For $x>0$, it is given by
\begin{equation}\label{b.vor}
D(x) =x\left(\log x
+2\gamma-1\right)+\frac{1}{4}+\sum_{n=1}^{\infty}d(n)\left(\df{x}{n}\right)^{1/2}M_1(4\pi\sqrt{nx}),
\end{equation}
where $\gamma$ denotes Euler's constant and $M_1(z)$ is defined by
\begin{equation}\label{b.c.9}
M_{\nu}(z) := -Y_{\nu}(z)-\df{2}{\pi}K_{\nu}(z).
\end{equation}
  Define the `error term' $\Delta(x)$ by
\begin{equation}\label{dirichlet}
D(x) =x\left(\log x
+2\gamma-1\right)+\frac{1}{4}+\Delta(x).
\end{equation}
Observe that $\Delta(x)$ is represented by the series of Bessel functions on the right-hand side of \eqref{b.vor}.  Finding the order of magnitude of $\Delta(x)$ as $x\to\infty$ is the famous, unsolved, notoriously difficult \emph{Dirichlet Divisor Problem}.  It is conjectured that for every $\epsilon>0$,
\begin{equation*}\label{delta}
\Delta(x)=O(x^{\frac14 +\epsilon}), \qquad x\to\infty.
\end{equation*}
Most attempts at finding upper bounds for $\Delta(x)$ as $x\to\infty$, in essence, depend on bounding partial sums of the infinite series in \eqref{b.vor}.  From \eqref{b5c} and \eqref{b6}, we see that the asymptotic formula \eqref{b5c} is a key ingredient in bounding $\Delta(x)$.  Currently, the best result is due to
M.~N.~Huxley \cite{huxley}, who proved that
$$\Delta(x)=O(x^{131/416+\epsilon}),$$
for each $\epsilon>0$.  Observe that $\frac{131}{416}=.3149\dots$. A survey of  much of the work that has been devoted to the \emph{Dirichlet Divisor Problem} can be found in \cite{circledivisorsurvey1}.

Other proofs of \eqref{b.vor} have been given by A.~Ivi\'{c} \cite{ivic1} and S.~Egger n\'{e} Endres and F.~ Steiner \cite{egg},

  A survey article on the life and work of Vorono\"{\dotlessi} has been written by J.~Steuding \cite{steuding}.

\section{An Identity of Ramanujan from His Lost Notebook}
 On a page published with his lost notebook
\cite[p.~335]{lnb}, Ramanujan recorded two identities involving
doubly infinite series of Bessel functions.  The first is associated with the classical \emph{Circle Problem} of Gauss, while the second is
connected with the \emph{Dirichlet Divisor Problem}. To state the latter formula as Ramanujan recorded it, first define
\begin{equation}\label{F}
F(x) = \begin{cases} [x], \quad &\text{if $x $ is not an
integer}, \\
x-\tfrac{1}{2}, \quad &\text{if $x$ is an integer},
\end{cases}
\end{equation}
where 
$[x]$ is the integral part of $x$.

 \begin{theorem}\label{lnbcosine}
\label{b.besselseries2} If $F(x)$ is defined by
\eqref{F}, then, for $x>0$ and $0<\theta<1$,
\begin{align}\label{b.c.11}
&\sum_{n=1}^{\infty}F\left(\frac{x}{n}\right)\cos(2\pi n\theta)=
\df{1}{4} -x\log(2\sin(\pi\theta))\notag \\
&+\df{1}{2}\sqrt{x}\sum_{m=1}^{\infty}\sum_{n=0}^{\infty}
\left(\df{M_1\left(4\pi\sqrt{m(n+\theta)x}\right)}{\sqrt{m(n+\theta)}}
+\df{M_1\left(4\pi\sqrt{m(n+1-\theta)x}\right)}{\sqrt{m(n+1-\theta)}}\right),
\end{align}
where $M_{1}(z) $ is defined by \eqref{b.c.9}.
\end{theorem}

For simplicity, assuming that $x$ is not an integer, by an elementary argument, we observe that
\begin{equation}\label{elemdiv}
D(x)=\sum_{n\leq{x}}\sum_{d|n}1=\sum_{dj\leq x}1=\sum_{d\leq x}\sum_{1\leq j\leq x/d}1=\sum_{d\leq{x}}\left[\df{x}{d}\right].
\end{equation}
   Setting $\theta=0$ in \eqref{b.c.11}, we see that the left side of \eqref{b.c.11} is a generalization of \eqref{elemdiv}, while the right-hand side contains a two-variable analogue of the series on the right side of \eqref{b.vor}.  We conjecture that Ramanujan established \eqref{b.c.11} in order to attack the \emph{Dirichlet Divisor Problem}, but he did not leave a record of any motivation or attempted proof.  The only proof of \eqref{b.c.11} has been given by Junxian Li, and the first and fourth authors, and  is difficult \cite{finalproblem}.

Extending the ideas of \cite{finalproblem}, the first author, Martino Fassina, Sun Kim, and the fourth author \cite{bfkz} established identities for certain finite sums of trigonometric functions in terms of Bessel functions.  We provide one example \cite[p.~6, Theorem 2.4]{bfkz}.
   
\begin{theorem}\label{even'}
Let $M_1(z)$ be defined by \eqref{b.c.9}. If $0<\theta,$ $\sigma <1$
and $x>0,$ then
\begin{align*}
&{\sum_{nm\leq x}}^{\prime}\cos(2\pi n\theta)\cos(2\pi m\sigma) \\
=&\frac14+\frac{\sqrt{x}}{4}\sum_{n,m\geq
0}\left\{\frac{M_1(4\pi\sqrt{(n+\theta)(m+\sigma)x})}{\sqrt{(n+\theta)(m+\sigma)}}
+\frac{M_1(4\pi\sqrt{(n+1-\theta)(m+\sigma)x})}{\sqrt{(n+1-\theta)(m+\sigma)}}\right.\nonumber
\\
&\qquad \qquad +\left.
\frac{M_1(4\pi\sqrt{(n+\theta)(m+1-\sigma)x})}{\sqrt{(n+\theta)(m+1-\sigma)}}+
\frac{M_1(4\pi\sqrt{(n+1-\theta)(m+1-\sigma)x})}{\sqrt{(n+1-\theta)(m+1-\sigma)}}\right\}.\nonumber
\end{align*}
\end{theorem}

\section{The Vorono\"{\dotlessi} Summation Formula}
In his two long papers \cite{voronoi}, Vorono\"{\dotlessi} established a more general result than \eqref{b.vor}, which is called the \emph{Vorono\"{\dotlessi} Summation Formula}.
If $f(x)$ is a suitable function, then Vorono\"{\dotlessi} proved that
\begin{align}\label{genvor}
\sum_{j=0}^{\infty}d(j)f(j)=&\int_{0}^{\infty}(2\gamma+\log x)f(x)dx \notag\\
&+\sum_{n=1}^{\infty}d(n)\int_{0}^{\infty}f(x)\left\{K_0(4\pi\sqrt{nx})-\frac12 \pi Y_0(4\pi\sqrt{nx})\right\}dx,
\end{align}
where $Y_0(z)$ and $K_0(z)$ are defined by \eqref{bessely-int} and \eqref{b66}, respectively, and $\gamma$ denotes Euler's constant. Alternative versions of \eqref{genvor} can be established in which the sums and integrals over $[0,\infty)$ are replaced by finite intervals. Vorono\"{\dotlessi}'s summation formula is reminiscent of the Poisson summation formula, and, indeed, proofs of \eqref{genvor} using the Poisson summation formula have been given; e.g., see \cite{landau} and \cite{hejhal}.
An especially short proof of \eqref{genvor} was given by N.~S.~Koshliakov \cite{kosh1} in 1929. Several applications of \eqref{genvor} were made in a series of papers by him.


We offer a special case of \eqref{genvor} in the spirit of \eqref{b.vor}.

K.~Soni \cite{soni}  and, subsequently, also S.~ Yakubovich \cite[p.~379]{yak1} showed that, for $y>0$,
\begin{gather}\sum_{n=1}^{\infty}\left\{K_0(2\sqrt{yn})-2y
    \df{\log(4\pi^2n/y)}{(4\pi^2n)^2-y^2}\right\}
    =\df18\log\left(\df{4\pi^2}{y}\right)-\df{\log\sqrt{y}}{y}-\df{y}{4}.\label{kab}
\end{gather}
However, this identity was first established by Vorono\"{\dotlessi} \cite[Equations (5), (6)]{voronoi}. It also appears in Ramanujan's Lost Notebook \cite[p.~254]{lnb}.

Another special case is given by Koshliakov's formula, given later in Theorem \ref{koshy}.

Define the generalized divisor function $d_k(n)$ by
$$\zeta^k(s)=\sum_{n=1}^{\infty}\df{d_k(n)}{n^s}, \quad \sigma >1.$$
In particular, M.~Jutila \cite{jutila}, A.~Ivi\'{c} \cite{ivic2}, and X.~Li \cite{li}  examined the ternary divisor problem, i.e., finding bounds for the `error term' for
$$\sum_{n\leq x}{d_3(n)}.$$
To study this problem, Jutila \cite{jutila} developed an analogue of \eqref{genvor}.


There now exist in the literature many analogues and generalizatioins of the Vorono\"{\dotlessi} summation formulas, established under a variety of conditions and in vastly different settings.  
 For example, see papers by E.~Assing \cite{assing}; E.~Assing and A.~Corbett \cite{assing-corbett};  D.~Banerjee, E.~M.~Baruch, and D.~Bump \cite{banerjee1};  D.~ Banerjee, E.~ M.~ Baruch, and E.~ Tenetov \cite{banerjee2}; B.~C.~Berndt \cite{V};  S.~Bettin and J.~B.~Conrey \cite{bettin-conrey}; A.~Corbett \cite{corbett};  A.~Dixit, B.~Maji and A.~Vatwani \cite{dmv};   D.~Goldfeld and X.~Li \cite{goldfeld3}, \cite{goldfeld1}, \cite{goldfeld2}; A.~Ichino and N.~Templier \cite{ichino}; E.~M.~Kıral and F.~Zhou \cite{kira};    N.~S.~Koshliakov \cite{koshwigleningrad}; A.~G.~F.~Laurin\u{c}ikas \cite{laur};   T.~Meurman \cite{meurman}; S.~D.~Miller and W.~Schmid \cite{miller1}, \cite{miller4}, \cite{miller3}, \cite{miller1}; S.~D.~Miller and F.~Zhou \cite{miller2};  Z.~ Qi \cite{qi}; N.~Templier \cite{templier}; C.~Nasim \cite{nasim}; T.~Watanabe \cite{watanabe}; T.~A.~Wong \cite{wong};   S.~Yakubovich \cite{yak3}, \cite{yak2}; V.~Zacharovas \cite{zak}; and F.~Zhou \cite{zhou}.   

\section{An Analogue of the Theta-Inversion Formula}
Consider a class of arithmetical functions studied by K.~Chandrasekharan and R.~Narasimhan \cite{cn1}.
Let $a(n)$ and $b(n)$, $1\leq n<\infty$, be two sequences of complex numbers, not identically 0.  Set
\begin{equation}\label{18}
\varphi(s):=\sum_{n=1}^{\infty}\df{a(n)}{\lambda_n^s}, \quad \sigma>\sigma_a; \qquad
\psi(s):=\sum_{n=1}^{\infty}\df{b(n)}{\mu_n^s}, \quad \sigma>\sigma_a^*,
\end{equation}
where $\sigma=\Re(s)$, $\{\l_n\}$ and $\{\mu_n\}$ are two sequences of positive numbers, each tending to $\infty$, and $\sigma_a$ and $\sigma_a^*$ are the (finite) abscissae of absolute convergence for $\varphi(s)$ and $\psi(s)$, respectively.  Assume  that $\varphi(s)$ and $\psi(s)$ have analytic continuations into the entire complex plane $\mathbb{C}$ and are analytic on $\mathbb{C}$ except for a finite set $\bf{S}$ of poles.
Suppose that for some $\delta>0$, $\varphi(s)$ and $\psi(s)$  satisfy a functional equation of the form
\begin{equation}\label{19}
\chi(s):=(2\pi)^{-s}\Gamma(s)\varphi(s)=(2\pi)^{s-\delta}\Gamma(\delta-s)\psi(\delta-s).
\end{equation}
Chandrasekharan and Narasimhan proved that the functional equation \eqref{19} is equivalent to the `modular' relation
\begin{equation}\label{modular}
\sum_{n=1}^{\infty}a(n)e^{-\lambda_n x}=\left(\df{2\pi}{x}\right)^{\delta}\sum_{n=1}^{\infty}b(n)e^{-4\pi^2\mu_n/x}+P(x), \qquad \Re(x)>0,
\end{equation}
where
\begin{equation*}
P(x):=\frac{1}{2\pi i}\int_{\mathcal{C}}(2\pi)^z\chi(z)x^{-z}dz,
\end{equation*}
where $\mathcal{C}$ is a curve or curves encircling all of $\bf{S}$.

 Chandrasekharan and Narasimhan also proved that \eqref{19} and \eqref{modular} are equivalent to the following identity for Riesz sums appearing on the left-hand side below.  Let $x>0$ and $\rho>2\sigma_a^*-\delta-\frac12$.  Then the functional equation \eqref{19} is equivalent to the \emph{Riesz sum} identity
\begin{gather}
\df{1}{\Gamma(\rho+1)}{\sum_{\lambda_n\leq x}}^{\prime}a(n)(x-\lambda_n)^{\rho}=
\left(\df{1}{2\pi}\right)^{\rho}
\sum_{n=1}^{\infty}b(n)\left(\df{x}{\mu_n}\right)^{(\delta+\rho)/2}J_{\delta+\rho}(4\pi\sqrt{\mu_n x})+Q_{\rho}(x),\label{21}
\end{gather}
where the Bessel function $J_{\nu}(z)$ is defined in \eqref{bb} and the prime $\prime$ on the summation sign on the left side indicates that if $\rho=0$ and $x\in\{\lambda_n\}$, then only $\tf12 a(x)$ is counted.  Furthermore,
 $Q_{\rho}(x)$ is  defined by
\begin{equation}\label{22}
Q_{\rho}(x):=\df{1}{2\pi i}\int_{\mathcal{C}}\df{\chi(z)(2\pi)^zx^{z+\rho}}{\Gamma(\rho+1+z)}dz,
\end{equation}
where $\mathcal{C}$ is a curve or curves encircling  $\bf{S}$.

In other words, if $a(n)$ and $b(n)$ are two sequences of arithmetical functions generated by Dirichlet series satisfying one of \eqref{19}, \eqref{modular}, or \eqref{21}, then the remaining two relations hold.

 The arithmetical functions $a(n)$ and $b(n)$ satisfying \eqref{19} and \eqref{modular} also satisfy a `modular' relation in which the exponential functions are, roughly, replaced by modified Bessel functions.  More precisely, the present authors proved the following theorem  \cite[Theorem 4.1]{BDGZ5}.

\begin{theorem}\label{maintheorem} Let $\Re(\nu)>-1$ and $\Re(c), \Re(r)>0$. Assume that the integral below converges absolutely.  Then,
\begin{gather}
\df{1}{2\pi r}\sum_{n=1}^{\infty}\df{a(n)}{(c^2+\lambda_n)^{(\nu-1)/2}}
K_{\nu-1}(4\pi r\sqrt{c^2+\lambda_n})\notag\\
=\df{1}{2\pi r^{\nu}c^{\nu-\delta-1}}
\sum_{n=1}^{\infty}\df{b(n)}{(r^2+\mu_n)^{(\delta-\nu+1)/2}}K_{\delta+1-\nu}(4\pi c\sqrt{r^2+\mu_n})\notag\\
+\int_0^{\infty}Q_{0}(x)(c^2+x)^{-\nu/2}K_{\nu}(4\pi r\sqrt{c^2+x})dx,\label{good}
\end{gather}
where $Q_0(x)$ is defined by \eqref{22}.
\end{theorem}

 The authors also established a slightly more general theorem, which we do not give here.  The first author's paper  contains the first statement and proof of Theorem \ref{maintheorem} \cite[pp.~343--344]{III}. One can also prove Theorem \ref{maintheorem} by using the Vorono\"{\dotlessi} summation formula \eqref{genvor} \cite{bcb}, \cite[p.~154]{V}.

  If we  let $\nu=1/2$ and appeal to the special case  
  \cite[pp.~79, 80]{watson}
\begin{equation}
K_{1/2}\left(z\right)=K_{-1/2}\left(z\right)=\sqrt{\frac{\pi}{2 z}}e^{-z},\label{I2}
\end{equation}
after some elementary manipulation, we obtain the identity
\begin{align*}\label{beauty}
\sum_{n=1}^{\infty}a(n)e^{-4\pi r \sqrt{c^2+\lambda_n}}=&
\,\,2^{3/2}c^{\delta+1/2}r\sum_{n=1}^{\infty}\df{b(n)}{(r^2+\mu_n)^{(\delta+1/2)/2}}K_{\delta+1/2}(4\pi c\sqrt{r^2+\mu_n})  \notag\\
&+2\pi r\int_0^{\infty}Q_0(x)\df{e^{-4\pi r\sqrt{c^2+x}}}{\sqrt{c^2+x}}dx,
\end{align*}
which might be compared with \eqref{modular}.

The beautiful asymmetry of the modular relation \eqref{modular} also characterizes \eqref{good} and other K-Bessel function identities in our survey.

We offer two examples to illustrate \eqref{good}.  Details for each example can be found in \cite{BDGZ5}.  Let $r_k(n)$ denote the number of representations of the positive integer $n$ as a sum of $k$ squares, where representations with different orders and different signs in the summands are regarded as distinct.  For example, if $k=2$,
$$ 5=(\pm 2)^2+ (\pm 1)^2=(\pm 1)^2+(\pm 2)^2,$$
and so $r_2(5)=8$.   In the notation of \eqref{18} and \eqref{19}, $\delta=k/2$ and $\lambda_n=\mu_n=n/2$. If we further define $r_k(0)=1$ and  replace $\nu$ by $\nu+1$, then \eqref{good} yields the identity
\begin{gather}
\sum_{n=0}^{\infty}\df{r_k(n)}{(c^2+n/2)^{\nu/2}}K_{\nu}(4\pi r\sqrt{c^2+n/2})\notag\\=
\df{1}{r^{\nu}c^{\nu-k/2}}\sum_{n=0}^{\infty}\df{r_k(n)}{(r^2+n/2)^{(k/2-\nu)/2}}K_{k/2-\nu}(4\pi c\sqrt{r^2+n/2}).
\label{watson4}
\end{gather}

The identity \eqref{watson4} was first proved by the first author, Y.~Lee, and J.~Sohn \cite[p.~39, Equation (5.5)]{bls}.  The special case $k=2$ was first proved by A.~L.~Dixon and W.~L.~Ferrar \cite[p.~53, Equation (4.13)]{dixonferrar} in 1934. Another proof for $k=2$ was given by Oberhettinger and Soni \cite[p.~24]{os}.

Second, recall that Ramanujan's arithmetical function $\tau(n)$ is defined by
\begin{equation}\label{taudelta}
\sum_{n=1}^{\infty}\tau(n)q^n:=\Delta(z)=q\prod_{n=1}^{\infty}(1-q^n)^{24},
\end{equation}
where $q=e^{2\pi i z}$ with $z\in \mathbb{H}$, that is $|q|<1.$
In the notation of \eqref{18} and \eqref{19}, $ \lambda_n=\mu_n=n$ and $\delta=12.$
Applying Theorem \ref{maintheorem} and replacing $\nu$ by $\nu+1$, we deduce that, for $\Re(\nu), \Re(c), \Re(r) >0$,
\begin{equation}\label{watson3}
\sum_{n=1}^{\infty}\df{\tau(n)}{(c^2+n)^{\nu/2}}K_{\nu}(4\pi r\sqrt{c^2+n})=
\df{1}{r^{\nu}c^{\nu-12}}\sum_{n=1}^{\infty}\df{\tau(n)}{(r^2+n)^{(12-\nu)/2}}K_{12-\nu}(4\pi c\sqrt{r^2+n}).
\end{equation}
Details for the proof of \eqref{watson3} can be found in either \cite{BDGZ5} or \cite{bls}, where the first proof of \eqref{watson3} was given.

Reference \cite{BDGZ5} also contains examples of \eqref{good} for $\sigma_k(n)$,  i.e.,
\begin{equation}\label{sumofdivisors}
\sigma_k(n):=\sum_{d|n}d^k;
\end{equation}
both even and odd primitive characters $\chi(n)$; and $F(n)$, the number of integral ideals of norm $n$ in an imaginary quadratic number field.

\section{Another Identity that Is Reminiscent of the Modular Relation \eqref{modular}}

Theorem \ref{thm3} below is an analogue of Theorem \ref{maintheorem}.  Theorem \ref{thm3} and a generalization of it were established in \cite[Theorems 4.1, 3.1]{proc}.

\begin{theorem}\label{thm3}  Assume that $\Re(\nu)>0$ and $\Re(s)>0$.  Also assume that  $\delta+\nu+1>\sigma_a^*>0$.  Suppose that the integral on the right side below converges absolutely.
Then,
\begin{align*}
\df{2}{s}\sum_{n=1}^{\infty}a(n)\lambda_n^{(\nu+1)/2}K_{\nu+1}(s\sqrt{\lambda_n})
=&2^{3\delta+\nu+1}s^{\nu}\pi^{\delta}\Gamma(\delta+\nu+1)
\sum_{n=1}^{\infty}\df{b(n)}{(s^2+16\pi^2\mu_n)^{\delta+\nu+1}}\notag\\
&+\int_0^{\infty}Q_{0}(x)x^{\nu/2}K_{\nu}(s\sqrt{x})dx,\label{30}
\end{align*}
where $Q_0(x)$ is defined by \eqref{22}.
\end{theorem}

  The first version of Theorem \ref{thm3} was established by the first author in \cite[p.~311]{paper1}; see also another paper by the first author \cite[p.~342]{III}.  Theorem \ref{thm3} was established via the  Vorono\"{\dotlessi} summation formula in \cite[p.~154]{V}.

Note that both Theorems \ref{maintheorem} and \ref{thm3} can be thought of as identities involving analogues of the Hurwitz zeta function $\sum_{n=0}^{\infty}(n+a)^{-s}$, where $\sigma >1$ and $a>0$.

In analogy with the examples in the previous section, we provide two corresponding illustrations for Theorem \ref{thm3}.  First, let $a(n)=r_k(n)$.  Then $\delta=k/2$ and $\lambda_n=\mu_n=n/2$. Set $s=2\pi\sqrt{2\beta}$, where $\Re(\sqrt{\beta})>0$. Then, for $\Re(\nu)>-1$,
\begin{equation}\label{39}
\sum_{n=0}^{\infty}r_k(n)n^{(\nu+1)/2}
K_{\nu+1}\left(2\pi\sqrt{n\beta}\right)=
\df{\beta^{(\nu+1)/2}\Gamma(\nu+1+k/2)}{2\pi^{k/2+\nu+1}}
\sum_{n=0}^{\infty}\df{r_k(n)}{(\beta+n)^{k/2+\nu+1}}.
\end{equation}
The identity \eqref{39} was first proved by A.~I.~Popov \cite[Equation (6)]{popov1935} in 1935.  See also \cite[p.~329, Corollary 4.6]{bdkz}, where \eqref{39} was established by a completely different method.

 Let $k=2$ and $\nu=-1/2$ above.  Revert to the notation in Theorem \ref{thm3} by setting $s=2\pi\sqrt{\beta}$. Then, for $\Re(s)>0$,
\begin{equation}\label{rh}
\sum_{n=0}^{\infty}r_2(n)e^{-s\sqrt{n}}=2\pi{s}\sum_{n=1}^{\infty}\df{r_2(n)}{(s^2+4\pi^2n)^{3/2}},
\end{equation}
in the notation of Hardy \cite{hardy}, who used \eqref{rh} to establish a lower bound
 for the error term (analogous to \eqref{dirichlet}) in Gauss's \emph{Circle Problem}.  (As indicated in \cite{hardy}, the identity \eqref{rh} also follows from an identity of Ramanujan.)

For another example of Theorem \ref{thm3}, set $a(n)=\tau(n)$.  We therefore can immediately deduce the identity
\begin{equation}\label{51}
\sum_{n=1}^{\infty}\tau(n)n^{(\nu+1)/2}K_{\nu+1}(s\sqrt{n})=
2^{36+\nu}s^{\nu+1}\pi^{12}\Gamma(13+\nu)\sum_{n=1}^{\infty}\df{\tau(n)}{(s^2+16\pi^2n)^{\nu+13}}.
\end{equation}
In particular, let $\nu=-1/2$.  Then, from
 \eqref{51} and \eqref{I2},
\begin{equation}\label{53}
\sum_{n=1}^{\infty}\tau(n)e^{-s\sqrt{n}}
=2^{36}\pi^{23/2}\Gamma\left(\frac{25}{2}\right)\sum_{n=1}^{\infty}\df{s\tau(n)}{(s^2+16\pi^2n)^{25/2}}.
\end{equation}
The identity \eqref{53} is originally due to Chandrasekharan and Narasimhan \cite[p.~16, Equation (56)]{cn1}.

\section{First Analogue of the Riesz Sum Identity \eqref{21}}

In the notation \eqref{18}, suppose that $\varphi(s)$ and $\psi(s)$ satisfy a functional equation of the form
\begin{equation}\label{IIIc}
\Gamma^2(\tfrac12 (s+1))\varphi(s)=\Gamma^2(\tfrac12 (2-s))\psi(1-s).
\end{equation}
Set, for arbitrary $\nu$,
\begin{equation}\label{IIId}
G_{\nu}(x)=Y_{\nu}(x)+\df{2}{\pi}e^{\pi i\nu}K_{\nu}(x).
\end{equation}
The following theorem is proved in \cite[Theorem 5.1, p.~333]{III}.

\begin{theorem} Let $\varphi(s)$ and $\psi(s)$ satisfy \eqref{IIIc}, and let $\rho$ denote a positive integer.  Then, for $\rho >2\sigma_a^*-\tfrac32$,
\begin{equation}\label{L}
\df{1}{\Gamma(\rho+1)}{\sum_{\lambda_n\leq x}}^{\prime}a(n)(x-\lambda_n)^{\rho}
=2^{-\rho}\sum_{n=1}^{\infty}b(n)\left(\df{x}{\mu_n}\right)^{(\rho+1)/2}G_{\rho+1}(4\sqrt{\mu_n x})
+Q_{\rho}(x),
\end{equation}
where $G_{\rho+1}(x)$ is defined by \eqref{IIId}, where $Q_{\rho}(x)$ is defined by
\begin{equation*}
Q_{\rho}(x):=\df{1}{2\pi i}\int_{\mathcal{C}}\df{\Gamma(z)\psi(z)x^{z+\rho}}{\Gamma(\rho+1+z)}dz,
\end{equation*}
and where $\mathcal{C}$ is a closed curve containing the integrand's poles on its interior.
\end{theorem}

We indicate one example, which we do not explicitly record. However, readers can readily substitute the relevant parameters into \eqref{L}.  Let $\chi$ denote a non-principal, odd, primitive, character modulo $k$.
Let $L(s,\chi)$ denote the associated Dirichlet $L$-function. Here,
$$\varphi(s)=\left(\df{\pi}{k}\right)^{-(s+1)}L^2(s,\chi).$$
Since $L(s,\chi)$ is an entire function, $Q_{\rho}(x)\equiv 0$. 
Then \cite[p.~71]{davenport}, 
$\lambda_n=\mu_m=\pi n/k$,
$$a(n)=\df{k}{\pi}\sum_{rs=k}\chi(r)\chi(s) \quad \text{and} \quad b(n)=\df{k}{\pi}\df{i\sqrt{k}}{G(\chi)}\sum_{rs=k}\overline{\chi}(r)\overline{\chi}(s),$$
and $G(\chi)$ is the Gauss sum
$$G(\chi)=\sum_{j=1}^k\chi(j)e^{2\pi i j/k}.$$

\section{A Second Analogue of the Riesz Sum Identity \eqref{21}}

In the notation \eqref{18}, suppose that $\varphi(s)$ and $\psi(s)$ satisfy a functional equation of the form
\begin{equation}\label{IIIa}
\Gamma(\tfrac12 s)\Gamma(\tfrac12(s-p))\varphi(s)=\Gamma(\tfrac12 (1-s))\Gamma(\tfrac12(p+1-s))\psi(p+1-s),
\end{equation}
where $p$ is an integer.  For arbitrary $\nu$ and any integer $p$, define
\begin{equation}\label{IIIb}
F_{\nu,p}(x):=\cos(\tfrac12(p+1)\pi)J_{p+\nu}(x)-\sin(\tfrac12(p+1)\pi)\left\{Y_{p+\nu}(x)-\frac{2}{\pi}e^{\pi i\nu}K_{p+\nu}(x)\right\}.
\end{equation}
In analogy with \eqref{21}, the following theorem holds \cite[Theorem 6.1, p.~336]{III}.

\begin{theorem} Let $\varphi(s)$ and $\psi(s)$ satisfy \eqref{IIIa}. Let $\rho$ denote a positive integer and $x>0$. Then, for $\rho>2\sigma_a^*-p-\tfrac32,$
\begin{align*}
\df{1}{\Gamma(\rho+1)}{\sum_{\lambda_n\leq x}}^{\prime}a(n)(x-\lambda_n)^{\rho}
=&2^{-\rho}\sum_{n=1}^{\infty}b(n)\left(\df{x}{\mu_n}\right)^{(\rho+p+1)/2}F_{\rho+1,p}(4\sqrt{\mu_n x})
+Q_{\rho}(x),
\end{align*}
where $F_{\rho+1,p}(x)$ is defined by \eqref{IIIb},  $Q_{\rho}(x)$ is defined by
\begin{equation*}
Q_{\rho}(x):=\df{1}{2\pi i}\int_{\mathcal{C}}\df{\Gamma(z)\psi(z)x^{z+\rho}}{\Gamma(\rho+1+z)}dz,
\end{equation*}
and $\mathcal{C}$ is a closed curve with the integrand's poles on its interior.
\end{theorem}

As in the previous section, for brevity, we provide one example without explicitly recording the identity.  Let
$$\varphi(s):=\pi^{-s}\zeta(s)\zeta(s-k), \quad k\neq 0,-1, \quad \sigma > 1, k+1, $$
where $\zeta(s)$ denotes the Riemann zeta function.  In the theorem above, let $p=k$.  Then the relevant parameters are
$$\lambda_n=\mu_n=\pi n \quad \text{and} \quad a(n)=b(n)=\sigma_k(n),$$
where $\sigma_k(n)$ is defined by \eqref{sumofdivisors}.

\section{Logarithmic Sums} In the notation of \eqref{18}, consider Dirichlet series $\varphi(s)$ and $\psi(s)$ satisfying a functional equation of the form
\begin{equation*}
\Gamma^m(s)\varphi(s)=\Gamma^m(\delta-s)\psi(\delta-s),
\end{equation*}
where $m$ is a positive integer and $\delta $ is a positive real number. (Note that for $m=1$ the notation above is different from that in \eqref{19}.)
 In \cite{IItams}, the first author established a general identity \cite[pp.~362, 365]{IItams} for
 \begin{equation}\label{IIa}
 S(a;q):=\df{1}{q!}{\sum_{\lambda_n\leq x}}^{\prime}a(n)\log^q(x/\lambda_n),
 \end{equation}
 where $q$ is a non-negative integer, and the prime $\prime$ on the summation sign on the left-hand side indicates that if $q=0$ and $\lambda_n=x$, then only $\tfrac12 a(x)$ is counted.  For $q=1$, as in \eqref{21}, ordinary Bessel functions arise in the general identity for \eqref{IIa}. We give one example.  Define
 \begin{equation*}
 \zeta_k(s):=\sum_{n=1}^{\infty}\df{r_k(n)}{n^s},\quad \sigma >1.
 \end{equation*}
 Then
 \begin{equation}\label{opp}
 \sum_{n\leq x}r_2(n)\log(x/n)=\pi x -\log x+\zeta_2 ^{\prime}(0)-\df{1}{\pi}\sum_{n=1}^{\infty}\df{r_2(n)}{n}J_0(2\pi\sqrt{nx}),
 \end{equation}
 where $\zeta_2 ^{\prime}(0)$ arises from the analytic continuation of $\zeta_2(s)$ into the entire complex plane. Oppenheim claimed that he had proved an identity for the left-hand side of \eqref{opp}, but he does not give it.

 For  $q=2$, the Bessel functions $K_{\nu}(x)$ and $Y_{\nu}(x)$ appear in the general identity for \eqref{IIa}.  Instead of offering the general identity, we provide two special cases, namely, for the divisor function $d(n)$ and $F(n)$, the number of integral ideals of norm $n$ in an algebraic number field. First, if $\gamma$ again denotes Euler's constant,
 \begin{align*}
 \sum_{n\leq x}d(n)\log(x/n) =&x(\log x-2+2\gamma)+\frac14\log(4\pi^2x)\\
 &+\df{1}{2\pi}\sum_{n=1}^{\infty}\df{d(n)}{n}\left\{Y_0(4\pi\sqrt{nx})+\df{2}{\pi}K_0(4\pi\sqrt{nx})\right\},
 \end{align*}
 which is originally due to Oppenheim \cite{oppenheim}; see also \cite[p.~371]{IItams}.  
 
 Second, let $h$ denote the class number of $K$, let $\Delta$ denote the discriminant of $K$, and let $\lambda=2\pi R/(w|\Delta|^{1/2})$, where $R$ is the regulator of $K$, and $w$ denotes the number of roots of unity in $K$.  Then,  
 \begin{align}\label{KK}
 \sum_{n\leq x}F(n)\log(x/n) =&\lambda hx+\zeta_K(0)\log x +{\zeta_K}^{\prime}(0)\notag\\
 &+\df{\sqrt{|\Delta|}}{2\pi}\sum_{n=1}^{\infty}\df{F(n)}{n}\left\{Y_0(4\pi\sqrt{nx/|\Delta|})+\df{2}{\pi}K_0(4\pi\sqrt{nx/|\Delta|})\right\},
 \end{align}
where $\zeta_K(s)$ is the Dedekind zeta function associated with the quadratic field $K$, that is,
 \begin{equation*}
 \zeta_K(s):=\sum_{n=1}^{\infty}\df{F(n)}{n^s},\quad \sigma >1,
 \end{equation*}
 which can be analytically continued to the entire complex plane.  See \cite[p.~374]{IItams} for a proof of \eqref{KK}.

 Returning to Theorem \ref{lnbcosine}, S.~Kim and the first author \cite[Theorem 5.1]{bk} proved a logarithmic analogue.  First, let 
 $$\zeta(s,a):=\sum_{n=0}^{\infty}(n+a)^{-s}, \quad \Re(s)>1, \,\, a>0,$$
 denote the Hurwitz zeta function.  Second, the generalized Stieltjes constants, $\gamma_n(a)$, $n\geq0$, are defined by the Laurent expansion
 \begin{equation*}
 \zeta(s,a)=\df{1}{s-1}+\sum_{n=0}^{\infty}\gamma_n(a)(s-1)^n.
 \end{equation*}
 In particular, $\gamma_0(1)=\gamma$, Euler's constant.  We now offer the aforementioned theorem from \cite{bk}.  
 
  \begin{theorem} Let $x>0$ and $0<\theta<1$.  Then
  \begin{align*}
  \sum_{n\leq x}&\log\left(\df{x}{n}\right)\sum_{r|n}\cos(2\pi r\theta)\\
  =&\df{\log(4\pi^2x)+\gamma}{4}-x\log(2\sin(\pi\theta))-\df18\left(\gamma_0(\theta)+\gamma_0(1-\theta)\right)\\
  &-\df{1}{4\pi}\sum_{n=0}^{\infty}\sum_{m=1}^{\infty}\left\{\df{M_0(4\pi\sqrt{m(n+\theta)x})}{m(n+\theta)}
  +\df{M_0(4\pi\sqrt{m(n+1-\theta)x})}{m(n+1-\theta)}\right\},
  \end{align*}
  where $M_0(z)$ is defined by \eqref{b.c.9}.
  \end{theorem}

\section{Theorem \ref{thm3} An Identity Involving Hypergeometric Functions} Recall that the rising or shifted factorial $(\alpha)_n$ is defined by
$$(\alpha)_n=\alpha(\alpha+1)\cdots(\alpha+n-1), \quad n\geq1; \qquad(\alpha)_0=1. $$
Let
$${_2F_1}(a,b;c;z):=\sum_{n=0}^{\infty}\df{(a)_n(b)_n}{(c)_nn!}z^n,\quad |z|<1,$$
denote the ordinary hypergeometric function.

In \cite{BDGZ5}, in the notation of \eqref{18} and \eqref{19}, the authors proved the following theorem.

\begin{theorem}\label{maintheorem2}
Assume that $\Re(\nu)>-1$,  $\Re(\sqrt{\alpha})>\Re(\sqrt{\beta})>0$, and  $\delta+\Re(\nu)+1>\sigma_a^*>0$. Let $Q_{\rho}(x)$ be defined by \eqref{22}, Bessel function $I_{\nu}(z)$, $K_{\nu}(z)$ be defined in \eqref{b1}, \eqref{b2}, and suppose that $Q_0(0)$ exists. Assume that the integral on the right side below converges absolutely.
Then,
\begin{align}
&\sum_{n=1}^{\infty}a(n)I_{\nu+1}\left(\pi \sqrt{\l_n}\left(\sqrt{\a}-\sqrt{\b}\right)\right)K_{\nu+1}\left(\pi \sqrt{\l_n}\left(\sqrt{\a}+\sqrt{\b}\right)\right)\nonumber\\
=&\frac{2(2\pi)^{-\delta}\Gamma(\nu+\delta+1)}
{\Gamma(\nu+2)}\sum_{n=1}^{\infty}\frac{b(n)}
{\sqrt{4\mu_n+\a}\sqrt{4\mu_n+\b}}
\left(\frac{\sqrt{4\mu_n+\a}-\sqrt{4\mu_n+\b}}{\sqrt{4\mu_n+\a}+\sqrt{4\mu_n+\b}} \right)^{\nu+1}\nonumber\\
&\times\left(\frac{1}{\sqrt{4\mu_n+\a}} +\frac{1}{\sqrt{4\mu_n+\b}}\right)^{2\delta-2}
{}_2F_{1}\left[\nu-\delta+2,1-\delta;\nu+2;
\left(\frac{\sqrt{4\mu_n+\a}-\sqrt{4\mu_n+\b}}{\sqrt{4\mu_n+\a}+\sqrt{4\mu_n+\b}} \right)^{2} \right]\nonumber\\
&+\frac{Q_{0}(0)}{2(\nu+1)}\left( \frac{\sqrt{\a}-\sqrt{\b}}{\sqrt{\a}+\sqrt{\b}}\right)^{\nu+1}+\int_{0}^{\infty}Q'_{0}(x)I_{\nu+1}\left(\pi \sqrt{x}\left(\sqrt{\a}-\sqrt{\b}\right)\right)K_{\nu+1}\left(\pi \sqrt{x}\left(\sqrt{\a}+\sqrt{\b}\right)\right)dx.\label{b13}
\end{align}
\end{theorem}

If we divide both sides of \eqref{b13} by $(\sqrt{\a}-\sqrt{\b})^{\nu+1}$, and let $\a\to \b$, then we obtain Theorem \ref{thm3}.  Details can be found in \cite{BDGZ5}.

Let $a(n)=r_k(n)$ in Theorem \ref{maintheorem}.  Then, after a considerable amount of work, we obtain the identity
\begin{align}\label{rkn1}
&\sum_{n=1}^{\infty}r_k(n)I_{\nu}\left(\pi \sqrt{n}\left(\sqrt{\a}-\sqrt{\b}\right)\right)K_{\nu}\left(\pi \sqrt{n}\left(\sqrt{\a}+\sqrt{\b}\right)\right)\nonumber\\
=&-\frac{1}{2\nu}\left( \frac{\sqrt{\a}-\sqrt{\b}}{\sqrt{\a}+\sqrt{\b}}\right)^{\nu}
+\frac{\Gamma(k/2+\nu)}{\pi^{k/2}2^{k-1}\Gamma(\nu+1)}
\sum_{n=0}^{\infty}\frac{r_k(n)}{\sqrt{n+\a}\sqrt{n+\b}}
\left(\frac{\sqrt{n+\a}-\sqrt{n+\b}}{\sqrt{n+\a}+\sqrt{n+\b}} \right)^{\nu}\nonumber\\
&\times\left(\frac{1}{\sqrt{n+\a}} +\frac{1}{\sqrt{n+\b}}\right)^{k-2}{}_2F_1\left[1-\frac{k}{2}+\nu,1-\frac{k}{2};\nu+1;
\left(\frac{\sqrt{n+\a}-\sqrt{n+\b}}{\sqrt{n+\a}+\sqrt{n+\b}} \right)^{2} \right].
\end{align}
The identity \eqref{rkn1} was first proved by Sun Kim and three of the present authors in \cite{bdkz}.

As a second example, let $a(n)=\tau(n)$ in
 Theorem \ref{maintheorem2}.  In this case, for $\Re(\nu)>-13/2$,
\begin{align}
&\sum_{n=1}^{\infty}\t(n)I_{\nu+1}\left(\pi \sqrt{n}\left(\sqrt{\a}-\sqrt{\b}\right)\right)K_{\nu+1}\left(\pi \sqrt{n}\left(\sqrt{\a}+\sqrt{\b}\right)\right)\nonumber\\
=&\frac{2(2\pi)^{-12}\Gamma(13+\nu)}{\Gamma(\nu+2)
}\sum_{n=1}^{\infty}\frac{\tau(n)}{\sqrt{4n+\a}\sqrt{4n+\b}}
\left(\frac{\sqrt{4n+\a}-\sqrt{4n+\b}}{\sqrt{4n+\a}+\sqrt{4n+\b}} \right)^{\nu+1}\nonumber\\
&\times\left(\frac{1}{\sqrt{4n+\a}} +\frac{1}{\sqrt{4n+\b}}\right)^{22}
{}_2F_{1}\left[\nu-10,-11;\nu+2;\left(\frac{\sqrt{4n+\a}-\sqrt{4n+\b}}{\sqrt{4n+\a}+\sqrt{4n+\b}} \right)^{2} \right].\label{c1}
\end{align}
Letting $\nu=-1/2$ in \eqref{c1}, we deduce that
\begin{gather*}
\sum_{n=1}^{\infty}\frac{\tau(n)}{\sqrt{n}}e^{-\pi\sqrt{n}(\sqrt{\a}+\sqrt{\b})}\sinh(\pi\sqrt{n}(\sqrt{\a}-\sqrt{\b}))\notag\\
=2\frac{3\cdot5\cdots21}{\pi^{11}}\sum_{n=1}^{\infty}\tau(n)
\left(\frac{1}{(4n+\b)^{23/2}}-\frac{1}{(4n+\a)^{23/2}}\right).
\end{gather*}

\section{Generalizations of Two Theorems of G.~N.~Watson}
The functional equation of the Riemann zeta function $\zeta(s)$ is given by \cite[p.~14]{edwards}
\begin{equation*}
\pi^{-s/2}\Gamma(s/2)\zeta(s)=\pi^{-(1-s)/2}\Gamma((1-s)/2)\zeta(1-s).
\end{equation*}
Replacing $s$ by $2s$, we see that it can be transformed into the form \eqref{19} with $\delta=1/2$, $a(n)=b(n)=1$, and $\lambda_n=\mu_n=n^2/2$.

As shown in \cite[Section 15]{BDGZ5}, the special case of Theorem \ref{maintheorem2} that arises from the Riemann zeta function is given by
\begin{align}\label{an1}
	&\frac{1}{4(\nu+1)}\left( \frac{\sqrt{\a}-\sqrt{\b}}{\sqrt{\a}+\sqrt{\b}}\right)^{\nu+1}+\sum_{n=1}^{\infty}I_{\nu+1}\left(\pi n\left(\sqrt{\a}-\sqrt{\b}\right)\right)K_{\nu+1}\left(\pi n\left(\sqrt{\a}+\sqrt{\b}\right)\right)\nonumber\\
	=&\frac{\Gamma(\nu+3/2)}{2\sqrt{2\pi}\Gamma(\nu+2)} \frac{\left(\sqrt{\a}-\sqrt{\b}\right)^{\nu+1}}{\left(\sqrt{\a}+\sqrt{\b}\right)^{\nu+2}}\cdot
	{}_2F_{1}\left(\nu+3/2,1/2;\nu+2;\left( \frac{\sqrt{\a}-\sqrt{\b}}{\sqrt{\a}+\sqrt{\b}}\right)^{2} \right)\nonumber\\
	&+\frac{\Gamma(\nu+3/2)}{\sqrt{\pi}\Gamma(\nu+2)}
	\sum_{n=1}^{\infty}\frac{\left(\sqrt{n^2+\a}
		-\sqrt{n^2+\b}\right)^{\nu+1}}{\left(\sqrt{n^2+\a}+\sqrt{n^2+\b}\right)^{\nu+2}}
	\cdot{}_2F_{1}\left[\nu+3/2,1/2;\nu+2;\left(\frac{\sqrt{n^2+\a}-\sqrt{n^2+\b}}{\sqrt{n^2+\a}+\sqrt{n^2+\b}} \right)^{2} \right],
\end{align}
where $\Re(\nu)>-1/2$ and $\Re(\sqrt{\alpha})>\Re(\sqrt{\beta})>0$.
Dividing both sides of \eqref{an1} by $(\sqrt{\a}-\sqrt{b})^{\nu+1}$, letting $\alpha\to\beta$, multiplying both sides of the resulting identity by $2\G(\nu+2)(2\sqrt{\beta})^{\nu+1}$, replacing $\nu$ by $\nu-1$ and $\beta$ by $z^2/(4\pi^2)$, and rearranging, for $\Re(z)>0$, we derive an intriguing identity of Watson \cite[Equation (4)]{watsonselfreciprocal}:
\begin{align*}\label{watson lemma}
\frac{1}{2}\G(\nu)+2\sum_{n=1}^\infty \left(\frac{1}{2}nz\right)^\nu K_\nu(nz)&=\sqrt{\pi}\G\left(\nu+\frac{1}{2}\right)z^{2\nu}\left\{\frac{1}{z^{2\nu+1}}+2
\sum_{n=1}^\infty \frac{1}{(z^2+4n^2\pi^2)^{\nu+\frac{1}{2}}}\right\}.
\end{align*}

Another consequence of \eqref{an1} is another identity of Watson \cite{watsonselfreciprocal}, \cite[Corollary 15.2]{BDGZ5}.  For
 $\Re( \beta)>0$,
\begin{equation*}
2\sum_{n=1}^{\infty}K_{0}(n\beta)=\pi\left\{\frac{1}{\beta}
+2\sum_{n=1}^{\infty}\left(\frac{1}{\sqrt{\beta^2+4\pi^2n^2}}
-\frac{1}{2n\pi}\right)\right\}+\gamma+\log\left(\frac{\beta}{2}\right)-\log 2\pi.
\end{equation*}

\section{Koshliakov's Formula and the Ramanujan--Guinand Formula}
In this section, we discuss related results of N.~S.~Koshliakov \cite{kosh1},  A.~P.~Guinand \cite{guinand}, and Ramanujan \cite{lnb}.  Guinand's Formula is probably the most famous of these theorems.  However, it is equivalent to a formula discovered by Ramanujan, which can be found in a partial manuscript published with his lost notebook \cite[pp.~253, 254]{lnb}.  We offer Ramanujan's formulation.

\begin{theorem} (\textbf{Guinand's Formula}) \cite[p.~253]{lnb}. Let $\sigma_k(n)$ be defined by \eqref{sumofdivisors}, let $\zeta(s)$
denote the Riemann zeta function, and let $K_{\nu}(z)$denote the modified Bessel function of order $\nu$ defined in \eqref{b2}.  If $\a$ and $\b$ are positive
numbers such that $\a\b=\pi^2$, and if $s$ is any complex number,
then
\begin{multline}\label{k.mainagain}
\sqrt{\a}\sum_{n=1}^{\infty}\sigma_{-s}(n)n^{s/2}K_{s/2}(2n\a)
-\sqrt{\b}\sum_{n=1}^{\infty}\sigma_{-s}(n)n^{s/2}K_{s/2}(2n\b)\\
=\df{1}{4}\Gamma\left(\df{s}{2}\right)\zeta(s)\{\b^{(1-s)/2}-\a^{(1-s)/2}\}
+\df{1}{4}\Gamma\left(-\df{s}{2}\right)\zeta(-s)\{\b^{(1+s)/2}-\a^{(1+s)/2}\}.
\end{multline}
\end{theorem}

The identity \eqref{k.mainagain} is equivalent to an identity
established by A.~P.~Guinand \cite{guinand} in 1955, about 35--40 years after Ramanujan discovered it.  Since \eqref{k.mainagain} is known as \emph{Guinand's Formula}, we have adhered to this convention.  It is closely allied to the Fourier expansion of nonholomorphic Eisenstein series on $\mathrm{SL}(2\mathbb{Z})$ or Maass wave forms \cite{maass}, and also to a famous result of  Selberg and Chowla on the Epstein zeta function \cite{sc}, which are discussed in Section \ref{selberg-chowla}.

A corollary of the Ramanujan--Guinand Formula  is a famous formula of Koshliakov \cite{kosh1} published in 1929.

\begin{theorem}\label{koshy} (\textbf{Koshliakov's Formula}) \label{kosh}
 Let  $d(n)$ denote the number of positive divisors
of the positive integer $n$. Then, if
 $\gamma$ denotes Euler's constant and  $a>0$,
\begin{gather}\label{k.kosh}
\gamma-\log\left(\df{4\pi}{a}\right)+4\sum_{n=1}^{\infty}d(n)K_0(2\pi
an) \notag\\ =\df{1}{a}\left(\gamma-\log(4\pi
a)+4\sum_{n=1}^{\infty}d(n)K_0\left(\df{2\pi n}{a}\right)\right).
\end{gather}
\end{theorem}

Koshliakov's proof depends upon the Vorono\"{\dotlessi} summation formula \eqref{genvor}.  However, Koshliakov's formula \eqref{k.kosh} is also originally due to Ramanujan in the same aforementioned two-page manuscript.
  It is clear from the manuscript that Ramanujan derived Koshliakov's formula from Guinand's formula, or more appropriately, the Ramanujan-Guinand formula. Later, K.~Soni \cite{soni}, C.~Nasim \cite{nasim}, and S. Yakubovich \cite[p.~374]{yak1} gave further proofs of \eqref{k.kosh}. Another representation for either side of \eqref{k.kosh} is given by an integral containing the Riemann $\Xi$-function in its integrand \cite[Eqeuation (1.12)]{dixit1}, whose equivalent form can actually be found in another paper of Koshliakov \cite[Equation (17)]{kosh1934}. The identity \eqref{k.kosh} should be compared with \eqref{kab}.

 Pages 253 and 254 of \cite{lnb} contain several other beautiful formulas.  We offer one of them as an example.

\begin{theorem}[p.~254]\label{k.tt}
If $a>0$,  $\gamma$ denotes Euler's constant, and $\sigma_k(n)$ is given by \eqref{sumofdivisors}, then
\begin{gather*}\label{k.3.20}
\mspace{-10mu} 2\sqrt a
\sum_{n=1}^\infty \sigma _{-1}(n) \sqrt{n} K_1(4\pi \sqrt{an}) \notag\\ =
-\dfrac{a^2}{2\pi} \sum_{n=1}^\infty
\dfrac{\sigma_{-1}(n)}{n(n+a)} + \dfrac{a}{2\pi}\left((\log a +
\gamma) \zeta(2) + \zeta^\prime(2)\right) + \dfrac{1}{4\pi} (\log
2a\pi +\gamma) + \dfrac{1}{48a\pi}.
\end{gather*}
\end{theorem}

Theorem \ref{k.tt} was only partially stated by Ramanujan.  The first author, Y.~Lee, and J.~Sohn \cite{bls} completed the formula and gave a proof for it.  All of Ramanujan's claims on these two pages can be found in the book by G.~E.~ Andrews and the first author \cite[Chapter 3]{bls}, as well as the aforementioned paper by  the first author, Lee, and Sohn \cite{bls}.  An equivalent formulation of Koshliakov's formula is due to the second author \cite{dixit1}. 

A generalization of the Ramanujan-Guinand formula was obtained by the second author, Kesarwani and Moll \cite[Theorems 1.4, 1.5]{dkmt}. It contains a new generalization of the modified Bessel function $K_z(x)$, see \cite[Equation (1.3)]{dkmt}. Several of its properties are established there. Yet another one was studied by Dixit and Kumar in \cite[Equation (1.1.37)]{super}. 

\section{The Chowla-Selberg Formula} \label{selberg-chowla} In 1949, S.~Chowla and A.~Selberg \cite{sc} offered what was later to be called the  \emph{Chowla-Selberg Formula}.  However, they did not provide a proof until 18 years later \cite{cs}.  To state their formula, first let $Q(m,n)=am^2+bmn+cn^2$ be a positive definite quadratic form with $a, b$, and $c$ real. The Epstein-zeta function $\mathbb{Z}(s,Q)$ is defined by
\begin{equation*}
\mathbb{Z}(s,Q):={\sum_{m,n=-\infty}^{\infty}}^{\prime}\{Q(m,n)\}^{-s},\quad \sigma >1,
\end{equation*}
where the prime $\prime$ suggests that the term corresponding to $(m, n)=(0, 0)$ should be omitted.

Then, the \emph{Chowla-Selberg Formula} is given by
\begin{align}\label{cs1}
a^{-s}\Gamma(s)\mathbb{Z}(s,Q)=&2\Gamma(s)\zeta(2s)+2k^{1-2s}\Gamma(s-\tf12)\zeta(2s-1)\notag\\
&+8k^{1/2-s}\pi^s\sum_{n=1}^{\infty}n^{s-1}\sigma_{1-2s}(n)\cos(n\pi b/a)K_{s-1/2}(2\pi kn).
\end{align}
Here,  $\sigma_k(n)$ is defined in \eqref{sumofdivisors}, $d=b^2-4ac$, $k=|d|/(4a^2)$, $\zeta(s)$ denotes the Riemann zeta-function, and $K_{\nu}(z)$ denotes the modified Bessel function, as usual.

Later proofs of \eqref{cs1} were provided by R.~A.~Rankin \cite{rankin}, P.~T.~Bateman and E.~Grosswald \cite{bateman-grosswald}, and Y.~Motohashi \cite{motohashi}.  The first author of the present paper established a generalization of \eqref{cs1} \cite{VI}.

\section{Real Analytic Eisenstein Series}\label{eisenstein}
Holomorphic Eisenstein series play a central role in the theory of modular forms. In particular, for the full modular group  $\Gamma=$ SL$_{2}(\mathbb{Z})$, the Eisenstein series for $k>2$ is defined by
\begin{align}\label{Eis}
E_{k}(z):=\sum_{\gamma\in \Gamma_{\infty}\backslash \Gamma}1|_k\gamma,
\end{align}
where
\begin{align*}\label{gammainf}
\Gamma_{\infty}=\left\{ \pm\begin{bmatrix}
1 & n \\
0 & 1
\end{bmatrix} \in \Gamma: n \in \mathbb{N}\right\},
\end{align*}
and the slash operator $|_k\gamma$ is defined as
\begin{align}
(f|_k\gamma)(z):=(cz+d)^{-k}f\left(\frac{az+b}{cz+d} \right), ~~~\gamma=\begin{bmatrix}
a & b \\
a & d
\end{bmatrix}\in  \textup{SL}_{2}(\mathbb{Z}), \quad z \in \mathbb{H}.\nonumber
\end{align}
Note that we can write \eqref{Eis} as a double sum, namely,
\begin{equation}\label{Eiss}
E_{k}(z)=\frac{1}{2}\sum_{\substack{-\infty< c,d<\infty\\(c,d)=1}}\frac{1}{(cz+d)^{k}}.
\end{equation}


For $\sigma>1$ and $z=x+iy \in \mathbb{H}$, a real analytic Eisenstein series on the full modular group SL$_{2}(\mathbb{Z})$ is defined by
 \begin{align}
E(s,z):=\frac{1}{2}\sum_{\gamma\in \Gamma_{\infty}\backslash \Gamma}\Im(\gamma z)^{s}.\nonumber
\end{align}
However, $E(s,z)$ is not a modular form on SL$_2(\mathbb{Z}).$
Similar to \eqref{Eiss}, $E(s,z)$ can also be expressed as a double sum
\begin{align}
E(s,z)=\frac{1}{2}\sum_{\substack{-\infty< c,d<\infty\\(c,d)=1}}\frac{y^s}{|cz+d|^{2s}}.\nonumber
\end{align}
Its modular properties and Fourier expansion, which we offer below, can be employed to prove that the Riemann zeta function does not vanish on the line Re$(s)=1$,  a key ingredient in one of  the proofs of the prime number theorem.

For any $\gamma=\begin{bmatrix}
a & b \\
c & d
\end{bmatrix} \in$ SL$_{2}(\mathbb{Z})$, we have the  modular relation \cite{zagier-eisenstein},
\begin{align*}
E(s,\gamma z)=E\left(s,\frac{az+b}{cz+d} \right)=E(s,z).
\end{align*}
The Fourier expansion of $E(s,z)$ is given in terms of a series involving $K_{\nu}(z),$ namely,
\begin{align}\label{functional}
E(s,z)=y^s+\frac{\psi(2s-1)}{\psi(2s)}y^{1-s}+\frac{4\sqrt{y}}{\psi(2s)}
\sum_{m=1}^{\infty}m^{s-\frac{1}{2}}\sigma_{1-2s}(m)K_{s-\frac{1}{2}}(2\pi my)\cos(2\pi m x ),
\end{align}
where 
\begin{align*}
\psi(s)=\pi^{-\frac{s}{2}}\Gamma\left( \frac{s}{2}\right)\zeta(s).
\end{align*}
Upon using  \eqref{b6}, we see that $E(s,z)$ has an analytic continuation to a meromorphic function on the entire complex plane in the variable $s.$ Also, if $\hat{E}(s,z)=\psi(2s)E(s,z)$, then the related functional equation $\hat{E}(s,z)=\hat{E}(1-s,z)$ holds. Indeed, note that
\begin{align}\label{Ecap}
\hat{E}(s,z)=\psi(2s)y^s+\psi(2s-1)y^{1-s}+4\sqrt{y}\sum_{m=1}^{\infty}m^{s-\frac{1}{2}}\sigma_{1-2s}(m)K_{s-\frac{1}{2}}(2\pi my)\cos(2\pi m x).
\end{align}
Employing the functional equations $\psi(s)=\psi(1-s)$, $\sigma_{\nu}(n)=n^{\nu}\sigma_{-\nu}(n)$, and $K_{\nu}(z)=K_{-\nu}(z)$ \cite[p.~79, Equation (8)]{watson} in \eqref{Ecap}, we obtain the functional equation
\begin{equation}\label{E}
\hat{E}(s,z)=\hat{E}(1-s,z).
\end{equation}
From the Fourier expansion \eqref{functional},  observe that $E(s,x+i y)$ is invariant under $x\to x+1$. Hence, we can write
\begin{align*}
E(s,z)=\sum_{m=-\infty}^{\infty}a_{m}(y)e^{2\pi i m x},
\end{align*}
where
\begin{align*}
a_{m}(y)=\left\{
	\begin{array}{ll}
		 y^s+\dfrac{\psi(2s-1)}{\psi(2s)}y^{1-s}, & \mbox{if } m = 0, \\
		2\sqrt{y}|m|^{s-\frac{1}{2}}\sigma_{1-2s}(|m|)K_{s-\frac{1}{2}}(2\pi |m|y),& \mbox{if } m \neq 0.
	\end{array}
\right.
\end{align*}

We provide a brief sketch of the well-known Kronecker limit formula, which is given by
\begin{align}\label{kroni}
    \left.\frac{\partial}{\partial s}E(s,z)\right|_{s=0}=\frac{1}{6}\log (y^6|\Delta(z)|,
\end{align}
where $\Delta(z)$ is defined in \eqref{taudelta}.
Recall \eqref{functional}, differentiate it with respect to $s$, and then employ \eqref{I2} in the series on the right-hand side to obtain
\begin{align}\label{szero1}
    \left.\frac{\partial}{\partial s}E(s,z)\right|_{s=0}=\log y-\frac{\pi}{3}y-\sum_{m=1}^{\infty}m^{-\frac{1}{2}}\sigma_{1}(m)e^{-2\pi my }\cos(2\pi m x).
\end{align}
On the other hand, if we use the definition of $\Delta(z)$, we can show that
\begin{align}\label{11111}
    \log(y^6|\Delta(z)|)&= 6\log y+\log(|\Delta(z)|)\nonumber\\
    &=6\log(y)-2\pi y-24\Re\left(\sum_{n=1}^{\infty}\sum_{m=1}^{\infty}\frac{q^{nm}}{m} \right)\nonumber\\
    &=6\log(y)-2\pi y-24\sum_{m=1}^{\infty}m^{-\frac{1}{2}}\sigma_{1}(m)e^{-2\pi my }\cos(2\pi m x).
\end{align}
Hence, from \eqref{szero1} and \eqref{11111}, we obtain \eqref{kroni}.

We conclude this section by giving an application of the functional equation \eqref{E} 
in the Rankin--Selberg method. This method is very helpful in obtaining the analytic continuation of $L$-functions associated with modular forms. We briefly explain the method   in obtaining an improved bound for Ramanujan's arithmetic function, $\tau(n)$, by employing the analytic continuation of  zeta functions attached to modular forms, in particular, to the cusp form $\Delta$-function, defined in \eqref{taudelta} (see \cite[Theorem 9.10]{kuro}).

To begin, let us first consider
\begin{align}\label{Ldelta}
L(\Delta,s):=\sum_{n=1}^{\infty}\frac{\tau(n)}{n^{s}}=\prod_{p}\left(1-\tau(p)p^{-s}+p^{11-2s} \right)^{-1}.
\end{align}
Note that the $L(s, \Delta)$ is absolutely convergent for $\sigma >13/2.$
Now, if we write \eqref{Ldelta} as
\begin{align*}
    L(\Delta,s)=\prod_{p}\left((1-\a_pp^{-s})(1-\b_pp^{-s})\right)^{-1},
\end{align*}
where $\a_p+\b_p=\tau(p)$ and $\alpha_p\beta_p=p^{11}$, then we can define,
\begin{align*}
    L(\Delta\otimes \Delta, s):=\prod_{p}\left((1-\a^2_pp^{-s})(1-\a_p\b_pp^{-s})^2(1-\b^2_pp^{-s})\right)^{-1}.
\end{align*}
Hence, $L(\Delta\otimes \Delta, s)=\zeta(2s-22)\sum_{n=1}^{\infty}\tau^2(n)n^{-s}$. Also,  $L(\Delta\otimes \Delta, s)$ has an analytic continuation to a meromorphic function on the entire complex plane, and
$$ \hat{L}(\Delta\otimes \Delta, s):=4(2\pi)^{11-2s}\Gamma(s)\Gamma(s-11) L(\Delta\otimes \Delta, s)$$
satisfies the functional equation $\hat{L}(\Delta\otimes \Delta, s)=\hat{L}(\Delta\otimes \Delta, s-11).$

In \cite{ran}, Rankin applied the ideas above to $\hat{L}(\Delta\otimes \Delta, s)$, and used its simple pole at $s=12$ to obtain the upper bound $\tau(n)=\mathcal{O}(n^{29/5})$, which improved the previous upper bound,  $\tau(n)=\mathcal{O}(n^6).$ We refer readers to \cite[Section 9.5]{kuro} for other applications in the theory of modular forms in which K-Bessel functions make an appearance.

\section{Applications in which the Bessel Functions $I_{\nu}(z)$ or $Y_{\nu}(z)$ appear}\label{ybessel}

Recall that the Bessel functions of imaginary argument $I_{\nu}(z)$  appear in Theorem \ref{maintheorem2}, while the Bessel functions $Y_{\nu}(z)$ defined in \eqref{b.Y} are focused in the formulas \eqref{b.vor} and \eqref{b.c.11}.  We complete our survey by mentioning two families of  examples in which each of these other Bessel functions appear, unaccompanied by their siblings.

First, as seen from their asymptotic expansions \eqref{b6} and \eqref{b5}, $K_{\nu}(z)$ and $I_{\nu}(z)$ play entirely different roles in analytic number theory.  In particular, the latter functions are prominent in describing the growth of partition functions.  Let $p(n)$ denote the number of partitions of the positive integer $n$.  For example, $p(5)=7$, because $5=4+1=3+2=3+1+1=2+2+1=2+1+1+1=1+1+1+1+1$. Improving on Hardy and Ramanujan's asymptotic formula for $p(n)$ \cite{circlemethod}, H.~Rademacher \cite{exact2} found an exact formula for $p(n)$, namely,
\begin{equation*}
p(n)=\df{2\pi}{(24n-1)^{3/4}}\sum_{k=1}^{\infty}\df{A_k(n)}{k}I_{3/2}\left(\df{\pi}{6k}\sqrt{24n-1}\right),
\end{equation*}
where $I_{3/2}$ is a modified Bessel function, $A_k(n)$ is the Kloosterman sum defined by
\begin{equation*}
A_k(n):=\sum_{\substack{0\leq h<k\\(h,k)=1}}e^{\pi i s(h,k)-2\pi inh/k},
\end{equation*}
and $s(h,k)$ denotes the Dedekind sum defined by
\begin{equation*}
s(h,k):=\sum_{r=0}^{k-1}\df{r}{k}\left(\df{hr}{k}-\left[\df{hr}{k}\right]-\df12\right).
\end{equation*}

There are several further asymptotic formulas for partition functions in which $I_{\nu}$ appears, for example, in recent papers by K.~Bringmann, B.~Kane, L.~Rolen and Z.~Tripp \cite{bkrt} and J.~ Iskander, V.~ Jain, and V.~Talvola \cite{ijt}.

	
Second,	Euler's famous formula linking  the Bernoulli numbers $B_{2n}, n\in\mathbb{N}$, with $\zeta(2n)$ is given by \cite[p.~5, Equation (1.14)]{temme-1996a}
	\begin{equation}\label{berexact}
		B_{2n}=\frac{2(-1)^{n+1}(2n)!\zeta(2n)}{(2\pi)^{2n}}.
	\end{equation}
In a beautiful paper \cite{zagier-1998a}, D.~Zagier obtained an analogue of this formula which features the modified Bernoulli numbers
\begin{equation}
B_{n}^{*} := \sum_{r=0}^{n} \binom{n+r}{2r} \frac{B_{r}}{n+r} \hspace{5mm} (n\in\mathbb{N}).
\nonumber
\end{equation}
He stumbled upon $B_n^*$ while devising a new proof of the Eichler--Selberg trace formula for the trace of Hecke operators $T_{\ell}$ acting
on modular forms on $\text{SL}_2{(\mathbb{Z})}$.
The numbers $B_{n}^{*}$ satisfy striking analogues of properties of $B_{n}$ as well. Zagier's analogue of \eqref{berexact} is given by \cite[Equation (10)]{zagier-1998a}
	\begin{align}
		B_{2n}^{*}   =  & -n+
		\sum_{\ell=1}^{\infty} \left( (-1)^{n} \pi Y_{2n}(4 \pi \ell) +
		\frac{1}{2 \sqrt{\ell}} \right) - \frac{1}{2} \zeta \left( \frac{1}{2} \right)
		\label{zagier-sum} \\
		 & + \sum_{\ell=1}^{\infty} \frac{1}{\sqrt{\ell(\ell+4)}}
		\left( \frac{\sqrt{\ell+4} - \sqrt{\ell}}{2} \right)^{4n},
		\nonumber
	\end{align}
where $Y_{n}(z)$ is defined by \eqref{bessely-int}.   The authors of \cite{dixit-2014a} studied properties of the Zagier polynomials defined by
\begin{equation*}
	B_{n}^{*}(x) = \sum_{r=0}^{n} \binom{n+r}{2r} \frac{B_{r}(x)}{n+r},\quad n>0,
\end{equation*}
where $B_{n}^{*}(0) = B_{n}^{*}$. These polynomials are analogues of the Bernoulli polynomials $B_n(x)$. It is well-known that $B_{2n}(x)$ and $B_{2n+1}(x)$, respectively, have the Fourier expansions \cite[p.~5]{zagier-1998a}
\begin{align}
	B_{2n}(x) &= 2(-1)^{n+1} (2n)! \sum_{m=1}^{\infty} \frac{\cos 2 \pi m x}
	{(2 \pi m)^{2n}},\quad 0\leq x\leq 1, n\geq 1,
	\label{ber-series1}\\
		B_{2n+1}(x) &= 2(-1)^{n+1} (2n+1)! \sum_{m=1}^{\infty} \frac{\sin 2 \pi m x}
	{(2 \pi m)^{2n+1}},	\label{ber-series2}
\end{align}
where the last formula holds for $0\leq x\leq 1$ when $n>0$, and for $0<x<1$ when $n=0$.

The authors of \cite{zagier3} studied analogues of \eqref{ber-series1} and \eqref{ber-series2} for the Zagier polynomials. We give below an analogue of \eqref{ber-series1} \cite[Theorem 1.2]{zagier3}.
\begin{theorem}
	\label{thm-main-0}
	Let $0 < x < 1$ and $n \in \mathbb{N} $. Define
	\begin{equation*}\label{gynx}
		g(y, r, x):=\frac{(y+1+x - \sqrt{(y-1+x)(y+3+x) } \,\, )^{2r}}{\sqrt{(y-1+x)(y+3+x)}}.
	\end{equation*}
	Let $Y_{n}(z)$ be defined in \eqref{bessely-int}, and denote by $U_{n}(x)$ the Chebyshev polynomial of the
	second kind, defined by
 \begin{align*}
U_{n}(x) = \frac{( x + \sqrt{x^{2}-1} )^{n+1} - ( x - \sqrt{x^{2}-1})^{n+1}}{2 \sqrt{x^{2}-1}}.
\end{align*}
 Then,
	\begin{eqnarray*}
		B_{2n}^{*}(x) &  = & (-1)^{n } \pi\sum_{m=1}^{\infty} Y_{2n}(4 \pi m) \cos(2 \pi m x) \\
		& & + \frac{1}{4} \left( U_{2n-1} \left( \frac{x+1}{2} \right) + U_{2n-1} \left( \frac{x}{2} \right) +
		U_{2n-1} \left( \frac{x-1}{2} \right) + U_{2n-1} \left( \frac{x-2}{2} \right)  \right) \\
		& & + \frac{1}{2^{2n+1}} \left(\sum_{m=1}^{\infty} g(m, n, x)+\sum_{m=1}^{\infty}g(m, n, 1-x)\right) .
	\end{eqnarray*}
\end{theorem}
There exists a similar formula for $B_{2n+1}^{*}(x)$ \cite[Theorem 1.3]{zagier3}. Letting $x\to1$ in Theorem \ref{thm-main-0} yields \eqref{zagier-sum}. This is analogous to obtaining \eqref{berexact} from \eqref{ber-series1} by substituting $x=1$ in the latter formula.

\section{Appreciation}

The authors thank University of Illinois Mathematics Librarians, Bernadette Braun and Becky Burner, for procuring several papers in our bibliography.

The first author is grateful to the Simons Foundation for its generous support. The second author sincerely thanks the Swarnajayanti Fellowship grant SB/SJF/2021-22/08 of SERB (Govt. of India) for its financial support.

\end{document}